\documentclass[12pt]{article}
\usepackage{amssymb,amsfonts}
\usepackage[english]{babel}
\usepackage{amsfonts,amsmath,amssymb,theorem}
\usepackage[unicode]{hyperref} % hyperlinks
\usepackage{color}
\usepackage{ifthen}
\usepackage{tikz}
\usetikzlibrary{calc}

\makeatletter
\def\@seccntformat#1{\csname the#1\endcsname.\ } % точка после номера раздела
\def\@biblabel#1{#1.} % формат номеров в списке литературы
\makeatother

\newenvironment{proof}[1][\hspace{-1.0ex}]%
{\par\addvspace{1mm}{\sc Proof\hspace{1.0ex}{#1}.} }%
{\quad$\blacktriangle$\par\addvspace{1mm}}

\newif\ifNoRemark
\def\addtheorem#1#2#3#4{
\ifthenelse{\equal{#2}{}}{}%
{\ifthenelse{\expandafter\isundefined\csname the#2\endcsname}{\newcounter{#2}}{}}
\newenvironment{#1}[1][\global\NoRemarktrue]% No Remark by default
{\par\addvspace{2mm plus 0.5mm minus 0.2mm}\noindent % Новый параграф без красной строки
{\bf #3}\ifthenelse{\equal{#2}{}}{}%
{\refstepcounter{#2}{\bf ~\csname the#2\endcsname}}%
{\bf \vphantom{##1}\ifNoRemark.\ \else\ (##1).\fi}\begingroup #4}%
{\endgroup\par\addvspace{1mm plus 0.5mm minus 0.2mm}\global\NoRemarkfalse}
\expandafter\newcommand\csname b#1\endcsname{\begin{#1}}
\expandafter\newcommand\csname e#1\endcsname{\end{#1}}
}

\addtheorem{theorem}{thrm}{Theorem}{\sl}
\addtheorem{lemma}{lmm}{Lemma}{\sl}
\addtheorem{proposition}{prps}{Proposition}{\sl}
\addtheorem{conjecture}{}{Conjecture}{\sl}
\addtheorem{corollary}{crll}{Corollary}{\sl}
\addtheorem{remark}{rmrk}{Remark}{\rm}

\title{On a Connection between the Switching
Separability of a Graph and That of Its Subgraphs}

\author{\href{http://arXiv.org/a/krotov_d_1}{Denis~S.~Krotov}%
\thanks{Sobolev Institute of Mathematics,
pr. Akad. Koptyuga 4, Novosibirsk, 630090 Russia}}

\date{}

\begin{document}
\maketitle
\begin{abstract}
A graph of order 
$n\geq 4$
is called {\it switching separable}
if its modulo-$2$ sum with some complete bipartite graph
on the same set of vertices
is divided into two mutually independent subgraphs,
each having at least two vertices.
We prove the following:
if removing any one or two vertices of a graph
always results in a switching separable subgraph,
then the graph itself is switching separable.
On the other hand, for every odd order greater than $4$,
there is a graph
that is not switching separable, but removing any vertex
always results in a switching separable subgraph.
We show a connection with similar facts on the separability
of Boolean functions and reducibility of $n$-ary quasigroups.

\end{abstract}

\setcounter{section}{-1}
\section{Introduction}

In the paper, we consider only simple graphs
(without loops and multiedges) and only induced subgraphs.
Let
$U$
be some set of vertices of a graph
$G=(V,E)$.
By a {\it switching}, or $U$-{\it switching}, of
$G$
one means the graph
$G_U = (V,E\bigtriangleup E_{U,V\setminus U})$,
where
$K_{U,V\setminus U}=(V,E_{U,V\setminus U})$
is the complete bipartite graph with parts
$U$, $V\setminus U$
(for generality reasons, we allow one of the parts to be empty).
It is easy to check that the relation
``$G'$ is a switching of $G$''
is an equivalence.
The set of switchings of some graph is called
a {\it switching class}.
It is known that there is a one-to-one correspondence
between the switching classes
and the so-called two-graphs \cite{krot6}.

We call a set
$W$
of vertices of a graph
$G=(V,E)$
{\it isolabde} if
$2\leq |W|\leq |V|-2$
and some switching of $G$ does not contain edges connecting
$W$
with
$V\setminus W$.
We call a graph of order
$n$
{\it switching separable}
if the vertex set includes an isolable subset.
Below, we will omit the word ``switching'' before ``separable.''

\begin{remark}
If a graph is separable, then all its switchings,
as well as its edge complement, are separable too.
All the graphs of order $4$ are separable.
\end{remark}

The goal of the current paper is to study relations
between the separability of a graph
and the separability of its subgraphs.
The research is motivated by connections of the topic with
the reducibility of $n$-ary quasigroups and the separability
of Boolean functions, which are discussed in Section~3.
In Sections~1 and~2, we will prove the following
two theorems:

\begin{theorem}\label{t1} {\it
Assume that all the subgraphs of orders
$n-1$ and $n-2$
of a graph
$G$
of order
$n$
are separable; then
the graph
$G$
is separable.
}
\end{theorem}

\begin{theorem}\label{t2} {\it
For every odd
$n\geq 5$,
there exists a non-separable graph of order
$n$
such that all its subgraphs of order
$n-1$
are separable.
}
\end{theorem}

The question about the even orders remains open,
but an exhaustive search has shown that there are
no similar examples of order $6$ or $8$.

\begin{conjecture}
From every non-separable graph of even order
one can obtain a non-separable subgraph by removing one vertex.
\end{conjecture}

The main result can be rephrased as follows:

\begin{corollary} {\it
From every non-separable graph,
one can obtain a non-separable subgraph by removing one or two vertices,
while removing one vertex is not sufficient in some cases.
}
\end{corollary}

The results were partially reported at the
IX International Workshop
``Discrete Mathematics and Its Applications,''
Dedicated to the 75th Anniversary of Academician O.~B.~Lupanov
(Moscow, June 18--23, 2007).

\section{Proof of Theorem \ref{t1}}\label{s:1}

Let
$\chi $
be the maximum order
of a non-separable proper subgraph of
$G$,
and let
$K$
be the vertex set of such a subgraph.
By the hypothesis,
$3\leq \chi\leq n-3$.

We first consider the case
$\chi>3$.
This inequality will be inexplicitly used
when establishing a contradiction with the non-separability of
$K$.
For four vertices
$a$, $b$, $c$, and~$d$
of
$G$,
denote by
$N(a,b;c,d)$
the number of edges of
$G$
among
$\{a,c\}$, $\{a,d\}$, $\{b,c\}$, and~$\{b,d\}$.

\begin{lemma}\label{l:1} {\it
A set
$W$
of vertices of
$G$
is isolable if and only if
for every distinct
$a,b$
from
$W$
and~$c,d$
from
$V\setminus W$
the number
$N(a,b;c,d)$
is even.
}
\end{lemma}
\begin{proof}
{\it Only if.}
Let
$W$
be an isolable set,
and let a set
$U$
define a corresponding separating switching
For every distinct
$a,b$
from
$W$
and~$c,d$
from
$V\setminus W$,
the graph $G$ contains exactly the same edges from
$\{a,c\}$, $\{a,d\}$, $\{b,c\}$, and~$\{b,d\}$
as the complete bipartite graph
$K_{U,V\setminus U}$.
It is easy to see that in a complete bipartite graph,
the number of edges that connect two pairs of vertices
is always even.

{\it If.}
Consider two nonadjacent vertices
$a$
from
$W$
and~$c$
from
$V\setminus W$
(if there are no such vertices, then
the $W$-switching isolates $W$).
Define the following four sets:
\begin{gather*}
W_0 = \{ b \in W \mid \{b,c\}    \in E \}, \qquad
V_0 = \{ d \in V\setminus W \mid \{a,d\}    \in E \},\\
W_1 = \{ b \in W \mid \{b,c\}\not\in E \}, \qquad
V_1 = \{ d \in V\setminus W \mid \{a,d\}\not\in E \}.
\end{gather*}
We claim that
$b$
from
$W_i$
and~$d$
from
$V_j$
are adjacent if and only if
$i+j=1$.
Indeed, if
$b=a$ or~$d=c$,
then the claim is straightforward from the definitions of
$W_i$ and~$V_j$;
otherwise, from the evenness of
$N(a,b;c,d)$.
So, taking
$U=W_0\cup V_0$,
we get
that the $U$-switching has no edges between
$W$ and~$V\setminus W$.
The proof of Lemma~\ref{l:1} is complete.
\end{proof}

\begin{lemma} \label{l:2}{ \it
For every five distinct vertices
$a$, $b$, $c$, $d$, and~$e$
of
$G$
the evenness of
$N(a,b;c,d)$ and~$N(a,b;c,e)$
implies the evenness of
$N(a,b;d,e)$.
}
\end{lemma}
\begin{proof}
We have
$$
N(a,b;d,e)=N(a,b;c,d)+N(a,b;c,e)-2|\{\{a,c\},
\{b,c\}\}\cap E|.
$$
%%% The proof of Lemma 2 is complete.
\end{proof}

Consider a vertex
$v\not\in K$.
By the definition of
$K$
the graph
$G|_{K\cup \{v\}}$
is separable;
i.e.,
$v$
belongs to some isolable set of vertices of
$G|_{K\cup \{v\}}$.
If this set has more than two vertices,
then the graph
$G|_{K}$
is also separable, which contradicts the definition of
$K$.
Hence, for every vertex
$v\not\in K$
there is a vertex
$u=u(v)\in K$
such that
$\{v,u\}$
is an isolable set in
$G|_{K\cup \{v\}}$.
Moreover,
$u(v)$
is defined uniquely
(if
$\{v,u\}$ and~$\{v,u'\}$
are isolable and
$u\neq u'$,
then from
Lemmas~\ref{l:1} and~\ref{l:2}
we get that
$\{v,u,u'\}$
is isolable,
which contradicts the non-separability of
$G|_{K}$).

\begin{proposition} {\it  \label{p:1}
For every two vertices
$v,v'\not\in K$
such that
$u(v)\neq u(v')$,
the number
$N(v,u(v);v',u(v'))$
is even.
}
\end{proposition}
\begin{proof}
As follows from theorem's hypothesis, the subgraph
$G|_{K\cup \{v,v'\}}$
is separable.
Let
$v$
belong to an isolable set $M$ of vertices
of this subgraph.
Consider subcases.

1.~
$|M\cap K|=0$, i.e., $M=\{v,v'\}$.
By Lemma~\ref{l:1}, for every different
$c,d\in K\setminus \{u(v)\}$,
the number
$N(v,v';c,d)$
is even.
The same is true for
$N(v,u(v);c,d)$,
because
$\{v,u(v)\}$
is isolable in
$G|_{K\cup \{v\}}$.
By Lemma~\ref{l:2},
$N(v',u(v);c,d)$
is also even, which implies that
$\{v,v',u(v)\}$
is isolable in
$G|_{K\cup \{v,v'\}}$
and 
$\{v',u(v)\}$
is isolable in
$G|_{K\cup \{v'\}}$.
The last contradicts to the uniqueness of
$u(v')$.

2.~
$|M\cap K|=1$.
If
$M=\{v,u(v)\}$,
then the claim of Proposition~\ref{p:1} follows from Lemma~\ref{l:1}.
Otherwise we have a contradiction to the uniqueness of
$u(v)$ or~$u(v')$.

3.~
$2\leq |M\cap K|\leq \chi-2$;
this contradicts the non-separability of
$G|_{K}$.

4.~
$|M\cap K|=\chi-1$.
The only vertex from
$K\setminus M$
is necessarily 
$u(v')$;
then the claim of
Proposition~\ref{p:1}
follows from
Lemma~\ref{l:1}.
%%% Proposition~1 is proved.
\end{proof}

\smallskip
Consider some vertex
$w$
from
$K$
with nonempty preimage

$u^{-1}(w)$.
Denote
$W = u^{-1}(w)\cup\{w\}$
and show that this set is isolable.
By Lemma~\ref{l:1},
this is equivalent to the evenness of
$N(a,b;c,d)$
for all distinct vertices
$a,b$
from
$W$
and
$c,d$
from
$V\setminus W$.
By Lemma~\ref{l:2},
it is sufficient to consider the case
$b=w=u(a)$.
If
$c$ and~$d$
belong to
$K$,
then the required evenness follows from the definition of
$u(a)$;
if
$d=u(c)$,
then, from
Proposition~\ref{p:1}.
The other cases are derived from these two using
Lemma~\ref{l:2}.
So, the separability of
$G$
in the case
$\chi>3$
is proved.

Consider the case
$\chi=3$.
Without loss of generality we may assume that
$G$
contains some isolated vertex
$o$
(otherwise, we can choose an arbitrary vertex as $o$
and make it isolated by switching the set of vertices adjacent to $o$).

\begin{proposition}\label{p:2} {\it
The graph
$G$
does not have a subgraph of type
$$
\left( \{o,a,b,c,d\},\left\{\{a,b\},\{b,c\},\{c,d\}\right\}\right).
$$
}
\end{proposition}

\begin{proof}
The non-separability of such a subgraph
contradicts theorem's hypothesis.
\end{proof}

\smallskip
We call two different vertices {\it twin}
if every other vertex
is adjacent to none or both of them.

\begin{proposition}\label{p:3} {\it
If
$v$ and $w$
are twin vertices, then 
$\{v,w\}$
is an isolable set.
}
\end{proposition}
\begin{proof}
Let
$U$
be the set of vertices adjacent to
$v$ and~$w$.
Then the $U$-switching of
$G$
does not contain edges connecting
$v$ or~$w$
with the other vertices.
\end{proof}

\smallskip
So, to establish the separability of
$G$,
it is sufficient to find twin vertices.
Consider a maximal sequence of vertices
$\bar u = (u_1, u_2, \ldots, u_t)$
that satisfies the following:

($*$)~
the vertices
$u_i$ and $u_j$, $1\leq i<j\leq t$,
are adjacent if and only if 
$i$
is odd.

If
$t=1$,
then the graph is empty and there is nothing to prove.
Suppose
$t>1$.
Let us show that
$u_{t-1}$ and~$u_t$
are twin vertices.
Seeking a contradiction, assume that this is not true.
Then there is a vertex
$w$
adjacent to exactly one of
$u_{t-1}$, $u_t$.
From
($*$)
we see that
$w$
is not in
$\bar u$.
We will show that the sequence
$\bar u$
is not maximal.

We first consider the case of odd
$t$.
Without loss of generality assume that
$w$
is adjacent with
$u_t$
(otherwise we interchange
$u_{t-1}$ and $u_t$, preserving the property ($*$)).
Note that
\begin{itemize}
\item
for every odd
$i$
less than 
$t$
the vertices
$u_i$ and~$w$
are adjacent;
otherwise the vertices
$o$, $u_{t-1}$, $u_i$, $u_t$, and~$w$
generate a forbidden subgraph
(Proposition~\ref{p:2});
\item
for every even
$i$
less than
$t-1$
the vertices
$u_i$ and~$w$
are not adjacent; otherwise the vertices
$o$, $u_i$, $w$, $u_{t-2}$, and~$u_{t-1}$
generate a forbidden subgraph
(Proposition~\ref{p:2}).
\end{itemize}
So, the sequence
$(u_1, u_2, \ldots, u_t,w)$
satisfies~($*$), which contradicts the maximality
of
$\bar u$
among the sequences satisfying~($*$).
The contradiction obtained proves that the vertices
$u_{t-1}$ and~$u_t$
are twin.

The case of an even
$t$
has a similar proof.
Theorem~\ref{t1} is proved.

%:333333333333333333333333333333333333333333333333333333333
\section{Proof of Theorem \ref{t2}}\label{s:2}

Let
$n\geq 5$
be odd.
Denote by
$G_n$
the graph with the vertex set
$V_n=\{v_i\}_{i=0}^{n-1}$
and the edges
$\{v_i,v_{i+j}\}$,
$j=1,\ldots,\left\lfloor  {n \over 4}\right\rfloor$
(here and below, the calculations with indexes are modulo
$n$).

\begin{proposition}\label{p:4} {\it
The graph
$G_n$
is not separable.
}
\end{proposition}
\begin{proof}
Denote
$m= \lfloor (n+1)/4 \rfloor$ and $u_i=v_{im}$.
Since
$n=4m\pm 1$,
the numbers
$m$ and~$n$
are relatively prime; hence
$\{u_i\}_{i=0}^{n-1} = V_n$.
Now consider an arbitrary subset
$A\subset V_n$
of cardinality at least $2$ and at most
$n-2$
and check that it is not isolable
(equivalently,
$V_n\setminus A$
is not isolable).
At least one of the following two cases takes place:

1.~
For some
$i$,
either
$u_i,u_{i+1}\in A$, $u_{i+2},u_{i+3}\not\in A$,
or
$u_i,u_{i+1}\not\in A$, $u_{i+2},u_{i+3}\in A$.
Then
$N(u_i,u_{i+1},u_{i+2},u_{i+3})=1$
(Fig.~1),
and,
by Lemma~\ref{l:1},
the set
$A$
is not isolable.

2.~
For some
$i$,
either
$u_i,u_{i+2}\in A$,
$u_{i+1}\not\in A$,
or
$u_i,u_{i+2}\not\in A$, $u_{i+1}\in A$.
Consider, for example,
the second subcase.
Note that every vertex
$u_j$
different from
$u_i$, $u_{i+1}$, $u_{i+2}$
is adjacent with exactly one of
$u_i$, $u_{i+2}$
(see Fig.~1).
Taking such
$u_j$
from
$A$,
we get that
$N(u_i,u_{i+2};u_{i+1},u_j)$
equals $1$ or $3$ (depending on
$n\bmod 4$).
Hence,
by Lemma~\ref{l:1},
the set
$A$
is not isolable.

Since every subset of vertices is not isolable,
the graph is not separable by the definition.
\end{proof}

%------------------------------  Figure 1 ------------
\begin{figure}[th]  %%% Figure:1
\def\iiii(#1-#2-#3-#4){\path[every node/.style={draw,rectangle,inner sep=1.2pt,fill=white}]
(#1) node{$\ a_i\ $}
(#2) node{$a_{\scriptscriptstyle i{+}1}$}
(#3) node{$a_{\scriptscriptstyle i{+}2}$}
(#4) node{$a_{\scriptscriptstyle i{+}3}$};}
\newcommand\dvsr{3}
\def\DefNode #1(#2){\node(#1) at ($(10,10)+sin(#2*360/\dvsr)*(10,0)+cos(#2*360/\dvsr)*(0,10)$){};}
\def\dfNM#1#2{\renewcommand\dvsr{#1}\node[rectangle,draw=white] at (10,10) {$\displaystyle{ n=#1 \atop m=#2}$};}
%55555555555555555555555555555555555555555555555555555555555
\mbox{}\hfill
\begin{tikzpicture}[scale=0.24, every node/.style={circle,draw,inner sep=2.0pt}]
\dfNM{13}{3}
\foreach \nam / \num in {a/0,b/1,c/2,d/3,e/4,f/5,g/6,h/7,i/8,j/9,k/10,l/11,m/12} {\DefNode \nam(\num)}
\draw[line width=0.7pt]
 (a)--(b)--(c)--(d)--(e)--(f)--(g)--(h)--(i)--(j)--(k)--(l)--(m)--(a)
    --(c)--(e)--(g)--(i)--(k)--(m)--(b)--(d)--(f)--(h)--(j)--(l)--(a)
    --(d)--(g)--(j)--(m)--(c) (l)--(b)--(e)--(h)--(k)--(a);
\draw[line width=2.2pt] (c)--(f)--(i)--(l);
\iiii(c-f-i-l)
\end{tikzpicture}
%66666666666666666666666666666666666666666666666666666666666
\ \ \ \ \ \ \ \ %
\begin{tikzpicture}[scale=0.24, every node/.style={circle,draw,inner sep=2.0pt}]
\dfNM{15}{4}
\foreach \nam / \num in {a/0,b/1,c/2,d/3,e/4,f/5,g/6,h/7,i/8,j/9,k/10,l/11,m/12,n/13,o/14} {\DefNode \nam(\num)}
\draw[line width=0.7pt]
 (a)--(b)--(c)--(d)--(e)--(f)--(g)--(h)--(i)--(j)--(k)--(l)--(m)--(n)--(o)--(a)
    --(c)--(e)--(g)--(i)--(k)--(m)--(o)--(b)--(d)--(f)--(h)--(j)--(l)--(n)--(a)
(a)--(d)--(g)--(j)--(m)--(a)(c)--(f)--(i)--(l)--(o)(b)--(e)--(h)--(k)--(n)--(b);
\draw[line width=2.2pt] (o)--(c);
\iiii(c-g-k-o)
\end{tikzpicture}
\hfill\mbox{}
%\onelinecaptionstrue \setcaptionmargin{30mm}
\caption{Examples for the graph
$G_n$,
the cases
$n\equiv 1\bmod 4$ and $n\equiv 3\bmod 4$
\label{fkr:1} }
\end{figure}
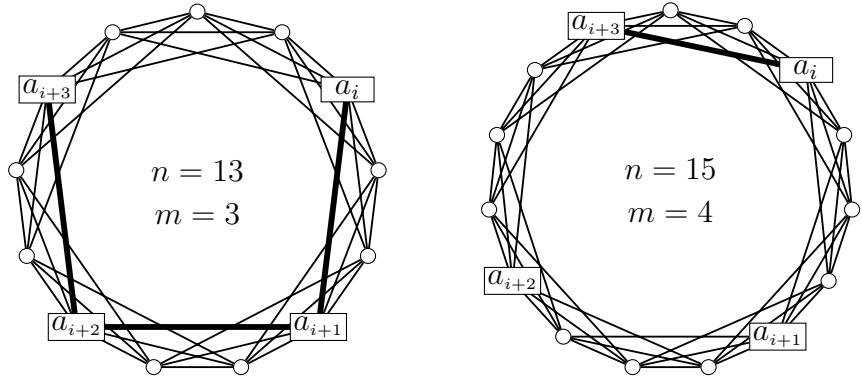
%------------------------------  Figure 1 ------------

\begin{proposition}\label{p:5} {\it
Removing any vertex in
$G_n$
results in a separable graph.
}
\end{proposition}
\begin{proof}
By symmetry, we may assume
that
$v_0$
is removed.
It is easy to see (Fig.~1)
that all other vertices except
$v_m$ and~$v_{-m}$
are divided into the vertices adjacent to
$v_m$
and the vertices adjacent to
$v_{-m}$.
As follows, the set
$\{v_m,v_{-m}\}$
is isolable (for the corresponding
$U$-switching, $U$ is defined as consisting
of
$v_m$ and all the vertices adjacent to $v_{-m}$).
%%% The proof of Proposition 5 is complete.
\end{proof}

\smallskip
So, Theorem~\ref{t2} is proved.
Note that, by Theorem~\ref{t1},
$G_n$
has a non-separable subgraph of order
$n-2$.
It can be shown similarly to the proof of
Proposition~\ref{p:4} that removing the vertices
$v_i$ and~$v_{i+m}$
gives such a subgraph.

%:444444444444444444444444444444444444444444444444444444444
\section{Graphs, boolean functions, quasigroups}\label{s:3}

In this section, we briefly discuss a connection of
the switching separability of graphs with similar properties
for Boolean functions and $n$-ary quasigroups.
Subgraphs of a graph correspond to so-called  retracts
of $n$-ary quasigroups and subfunctions of Boolean functions,
which are obtained by fixing some arguments.
In terms of retracts and subfunctions,
for $n$-ary quasigroups and Boolean functions
there hold theorems similar to
Theorems~\ref{t1} and~\ref{t2}
(the last is known for $n$-ary quasigroups only in the case
when the order is divisible by $4$).
Moreover, taking into account that,
using quadratic polynomials, one can map the graphs
into Boolean functions, and then to $n$-ary quasigroups of order $4$,
Theorem~\ref{t1}
is in fact a corollary of the corresponding theorem for quasigroups
\cite{krot2, krot4},
and Theorem~\ref{t2}, conversely, provides existence
of similar examples for
$n$-ary quasigroups of order
$4$~\cite{krot1}.

\subsection{Extended Boolean functions}

By an {\it extended Boolean function} we will mean
a partial Boolean function which is defined on the
binary $n$-tuples with even number of ones.
Note that an extended Boolean function can be treated
as usual Boolean function in $n-1$ arguments.
We call an extended Boolean function
$f$
in
$n$
arguments {\it separable}
if it can be represented as the sum of two Boolean functions
$f'$ and~$f''$
in
$n-2$
or smaller numbers of arguments,
the sets of arguments of
$f'$ and~$f''$
being disjoint (the bound
$n-2$
is rather natural:
in this case, to describe the function,
one have to define a smaller number
of point values of two decomposition functions
than
when listing the point values of the function itself).
By the {\it degree} of an extended Boolean function we mean
the minimal degree of a polynomial (over the field GF(2))
that represents this function.
By ``quadratic'' we mean ``of degree at most two.''
By the {\it graph of a quadratic polynomial} we mean
the graph on the set of arguments
such that two vertices are adjacent if and only if
the polynomial includes the product of the corresponding variables.

\begin{lemma}  \label{l:3} {\it
The set of graphs that correspond to the representations
of an extended Boolean function by quadratic polynomials
forms a switching class.
}
\end{lemma}
\begin{proof}
Every (in particular, quadratic) extended Boolean function
$f$
in
$n$
arguments,
being a Boolean function in the first
$n-1$
of its arguments,
is uniquely represented as
$$
f(x_1,\ldots,x_{n-1},x_{n}) = p(x_1,\ldots,x_{n-1}),
$$
where
$p$
is a polynomial.

Every polynomial
$r$
in
$n$
variables
$x_1,\ldots,x_{n-1},x_n$
can uniquely be represented as
$$
 q(x_1,\ldots,x_{n-1})
 + (x_1 + \dots + x_{n-1} + x_n) l(x_1,\ldots,x_{n-1}),
 $$
where
$q$ and $l$
are polynomials in
$x_1,\ldots,x_{n-1}$;
moreover, if
$r$
is quadratic, then
$q$
is quadratic and
$l$
is linear.
Since
$x_1+\dots+x_{n-1}+x_n=0$
over the domain of definition of extended Boolean functions,
the polynomial
$q$
is the same for all polynomials
that represent the same extended Boolean function.
It is easy to verify that the addition of
$(x_1+\dots+x_{n-1}+x_n) l(x_1,\ldots,x_{n-1})$
with a linear
$l$
results in a switching of the corresponding graph,
precisely, in the $U$-switching,
where
$U$
is the set of variables that $l$ essentially depends in.
Now the claim of Lemma~\ref{l:3} is straightforward.
\end{proof}

\begin{lemma}\label{l:4} {\it
A quadratic extended Boolean function is separable
if and only if the graphs of the quadratic polynomials
representing this function
are separable.
}
\end{lemma}
\begin{proof}
Straightforwardly, the separability of the graph
yields the separability of the function
represented by the polynomial with this graph.

To prove the converse, taking into account the previous lemma,
it is sufficient to prove that
for a separable quadratic extended Boolean function
$f$,
the elements
$f'$ and~$f''$
of some of its decompositions are also quadratic.
Let
$$
f(x_1,\ldots,x_n) = f'(\bar y)+f''(\bar z),
$$
where
$\bar y$ and $\bar z$
are disjoint collections of variables from
$x_1,\ldots,x_n$.
Consider the representation of
$f'(\bar y)+f''(\bar z)$,
as a total Boolean function, in the form
$$
f'(\bar y)+f''(\bar z)
= q(x_1,\ldots,x_{n-1})
+ (x_1+\dots+x_{n})l(x_1,\ldots,x_{n-1}).
$$
Since
$f$
is quadratic,
$q$
is quadratic too (see the proof of Lemma~\ref{l:3}).
Divide the polynomial
$l$
into the sum of two polynomials
$l_1+l_2$,
where
$l_1$
is linear and
$l_2$
consists of monomials of degree at least~$2$.
We have
$$
f'(\bar y) + f''(\bar z) = q(x_1,\ldots,x_{n-1})
+ \sum_{i=1}^{n} x_i l_1(x_1,\ldots,x_{n-1})
+ \sum_{i=1}^{n} x_i l_2(x_1,\ldots,x_{n-1}).
$$
We see that the last summand consists of
monomials of degree at least $3$
(indeed, if $l_2$ contains a monomial
$x_i x_j$,
then its product with
$x_i+x_j$
is zero, while the product with the other variables
consists of monomials of degree $3$).
So, removing this summand is equivalent to removing
the monomials of degree more than $2$ in the polynomial
representation of
$f'$ and~$f''$,
which results in the identity
$$
g'(\bar y) + g''(\bar z) = q(x_1,\ldots,x_{n-1})
+ \sum_{i=1}^n x_i l_1(x_1,\ldots,x_{n-1})
$$
for some quadratic functions
$g'$ and~$g''$.
Obviously,
$g'(\bar y)+g''(\bar z)$
also coincides with
$f$
on the domain of definition of an extended Boolean function.
So,
we get a quadratic representation of
$f$
that corresponds to a separable graph.
By Lemma~\ref{l:3},
all other quadratic representations also correspond
to separable graphs.
%%% This completes the proof of Lemma 4.
\end{proof}

\subsection{$n$-Ary Quasigroups}

Let
$\Sigma$
be a nonempty set.
An $n$-ary operation
$Q : \Sigma^n \to \Sigma$
is called an {\it $n$-ary quasigroup} of order
$|\Sigma|$,
if in the equation
$x_0=Q(x_1,\ldots,x_n)$
the values of any
$n$
variables uniquely determine the value of the remaining variable.
(Strictly speaking, an $n$-ary quasigroup is the pair
$(\Sigma,Q)$,
but we use a standard simplification of the terminology.)
As follows from the definition,
an $n$-ary quasigroup is invertible in every argument;
in the case of a finite order this property can be taken
as a definition.

We will use the following predicate notation for an $n$-ary
quasigroup:
$$Q\langle x_0,x_1,\ldots,x_n\rangle ~
\Leftrightarrow ~ x_0=Q(x_1,\ldots,x_n).$$
Often, the predicate notation occurs more convenient
than the functional one because of the symmetry with respect to
all $n+1$ variables.
If one fix the values of some
$m\in\{1,\ldots,n\}$
arguments in the predicate
$Q\langle \ldots\rangle$,
then the
$(n{+}1{-}m)$-ary
predicate obtained corresponds to
an $(n{-}m)$-ary quasigroup,
which is called a {\it retract} of
$Q$.
An $n$-ary quasigroup is called {\it permutably reducible}
(below, simply {\it reducible})
if it can be represented as a repetition-free composition
of two quasigroups of smaller arity,
where the order of variables in the composition
may differ from their original order.

\begin{remark}
In the literature,
by the reducibility (without ``permutable'')
one often means
that the quasigroup is decomposable into a composition
with the same order of variables.
In the Russian-language papers,
including the original of this contribution,
the permutably reducible quasigroups
are also known as {\em separable}.
\end{remark}

Let
$\Sigma=\{[0,0],[0,1],[1,0],[1,1]\}$
be the set of binary pairs,
and let
$\lambda$
be some extended Boolean function in
$n+1$
arguments.
The predicate
$$
Q_{\lambda} \Bigl\langle [x_0,y_0],\ldots,[x_n,y_n] \Bigr\rangle
~~ \Leftrightarrow ~~
\left\{ { |x_0 + \dots + x_n| = 0,
\hfill\mbox{}
\atop |y_0 + \dots + y_n| = \lambda(x_0,\dots,x_n) }\right.
$$
corresponds to an $n$-ary quasigroup
$Q_{\lambda}$
(addition is performed modulo $2$);
this construction is a partial case of
so-called wreath product of $n$-ary quasigroups,
in our case, of trivial $n$-ary quasigroups of order $2$.
(The term ``wreath product'' for $n$-ary quasigroups
does not agree with the known wreath product for groups;
so one should care to avoid misunderstanding using this term.)

\begin{lemma}[\cite{krot1}] \label{l:5} {\it
The reducibility of the $n$-ary quasigroup
$Q_{\lambda}$
is equivalent to the separability of the extended Boolean function
$\lambda$.
}
\end{lemma}

So, the separability of graphs is close related with
the reducibility of $n$-ary quasigroups,
at least within the framework of the following construction:
starting from a graph of order $n+1$,
we construct a quadratic extended Boolean function
$\lambda$ (the edges correspond to monomials of degree $2$;
the linear part is chosen arbitrarily);
then, we construct the $n$-ary quasigroup
$Q_\lambda$
of order~$4$;
furthermore,
we can construct an $n$-ary quasigroup of order
$4k$
for every % ¤«п «оЎ®Ј®
$k$,
including infinite,
using the direct product with the $n$-ary quasigroup
$P(x_1,\ldots,x_n)=x_1\star\ldots\star x_n$,
where % Ј¤Ґ
$\star$
is a commutative group operation.

In this sequence, the separability (reducibility)
of every element is equivalent to the  separability (reducibility)
of all other elements.
Moreover, the separability of a subgraph
is equivalent to the separability of the corresponding
subfunctions of the extended Boolean function
and the reducibility of the corresponding retracts
of the $n$-ary quasigroup.
This means that
Theorem~\ref{t1}
is a corollary of
Theorem~\ref{t3} below,
while
Theorem~\ref{t4}
follows from
Theorem~\ref{t2}.
In support of Sections~\ref{s:1}  and~\ref{s:2}
we note that the arguments there are essentially
easer and more readable than
the proofs for quasigroups.

In conclusion, we quote the known theorems for
$n$-ary quasigroups related to Theorems~\ref{t1} and~\ref{t2}
considered in the current paper.
For an $n$-ary quasigroup
$Q$,
denote by
$\chi(Q)$
the maximum arity of its irreducible retract.

\begin{theorem}\label{t3} {\it
If
$\chi(Q)<n-2$,
then the $n$-ary quasigroup
$Q$
is reducible.
If
$\chi(Q)=n-2$
and the order of
$Q$
is a prime integer, then
$Q$
is reducible.
}
\end{theorem}

The case
$2 < \chi(Q) < n-2$
was considered in
\cite{krot2};
the case
$\chi(Q)=2$
for order
$4$,
in \cite{krot4};
the general case with
$\chi(Q)=2$
and the case
$\chi(Q)=n-2$
for prime orders was proved in
\cite{krot5}.
Theorem~\ref{t3}
is useful for
inductive characterization of classes of
$n$-ary quasigroups;
for example,
in the proof that every $n$-ary quasigroup of order $4$
is semilinear (that is, equivalent to some
$Q_\lambda$)
or reducible \cite{krot4}, it used that a minimal
counterexample must have an irreducible
semilinear
$(n-1)$-ary
or
$(n-2)$-ary retract.

\begin{theorem}\label{t4} {\it
For every even
$n\geq 4$
and every
$k\geq 1$,
there exists an irreducible $n$-ary quasigroup
$Q$
of order
$4k$
such that
$\chi(Q)=n-2$ {\rm \cite{krot1}}.
For every 
$n\geq 3$ and $k\geq 4$,
there exists an irreducible
$n$-ary quasigroup of order
$k$
with
$\chi(Q)=n-1$ {\rm \cite{krot3}}.
}
\end{theorem}

The following cases remain uninvestigated
from the point of view of existence of
irreducible $n$-ary quasigroups:
\begin{itemize}
\item
$\chi(Q)=n-2$,
odd
$n$
(this case is connected with the conjecture formulated in the introduction),
non-prime orders;
\item
$\chi(Q)=n-2$,
arbitrary $n\geq 4$,
non-prime orders that are not divisible by $4$;
\end{itemize}

There is an example of irreducible $4$-ary quasigroup $Q$ of order $6$
with
$\chi(Q)=2$.
\section*{Acknowledgments}

The author is grateful to
A.~N.~Glebov, V.~N.~Potapov, A.~V.~Pyatkin,
and the anonymous referee
for their interest to this work and useful comments,
in particular, finding a gap in the proof of
Theorem~\ref{t1} in a preliminary version of this manuscript.

%The author was supported by
%the Federal Target Grant
%``Scientific and Educational Personnel
%of Innovation Russia'' for 2009--2013
%(government contract no.~02.740.11.0429)
%and
%the Russian Foundation for Basic Research
%(projects nos.~10--01--00616 and~10--01--00424).

\end{document}

\title{О связи свитчинговой разделимости графа и его подграфов%
%\renewcommand\thefootnote{}
%\thanks{Работа выполнена при поддержке ФЦП «Научные и научно"=педагогические кадры инновационной России» на 2009-2013 гг. (гос. контракт № 02.740.11.0429) и Российского фонда фундаментальных исследований (проект 08-01-00673)
%}
\thanks{Результаты работы частично докладывались на IX Международном семинаре <<Дискретная математика и ее приложения>>,
посвященном 75-летию со дня рождения академика О.\,Б.\,Лупанова (Москва, 18-23 июня 2007\,г.)}
}
\def\today{}
\author{\href{http://arXiv.org/a/krotov_d_1}{Д.~С.~КРОТОВ}
\thanks{Адрес автора: Институт математики им. С. Л. Соболева СО РАН,
проспект Академика Коптюга 4, Новосибирск, 630090, Россия
(e-mail: krotov@math.nsc.ru)}
 }
\maketitle

УДК: 519.173.1

\begin{abstract}
Граф порядка $n\geq 4$ называется свитчингово разделимым, если его сумма по модулю два с некоторым полным двудольным графом на том же множестве вершин разделена на два несвязанных между собой подграфа на двух или более вершинах. Доказано, что если удалением одной или двух вершин из данного графа мы получаем только свитчингово разделимые подграфы, то и сам граф свитчингово разделим. С другой стороны, существует граф любого нечетного порядка, который сам не является свитчингово разделимым, а удаление любой вершины приводит к свитчингово разделимому подграфу. Показана связь с аналогичными фактами для разделимости булевых функций и $n$-арных квазигрупп.
\end{abstract}
%:111111111111111111111111111111111111111111111111111111111
\setcounter{section}{-1}
\section{Введение}
В данной работе рассматриваются только простые графы (без кратных дуг и без петель),
и только индуцированные подграфы.
Пусть $U$ --- некоторое множество вершин графа $G=(V,E)$.
\emph{Свитчингом}, или \emph{$U$-свитчингом}, графа $G$ называется граф
$G_U = (V,E\bigtriangleup E_{U,V\setminus U})$, где $K_{U,V\setminus U}=(V,E_{U,V\setminus U})$
есть полный двудольный граф с долями $U$, $V\setminus U$ (для общности будем считать, что одна из долей может быть пустой).
Легко убедиться, что отношение <<$G'$ есть свитчинг $G$>> является эквивалентностью. Множество свитчингов одного графа называется \emph{свитчинговым классом}. Известно взаимнооднозначное соответствие между свитчинговыми классами и так называемыми два-графами \cite{two}.

Множество $W$ вершин графа $G=(V,E)$ назовем \emph{отделимым}, если $2\leq |W|\leq |V|-2$ и некоторый свитчинг графа не содержит ребер соединяющих $W$ с $V\setminus W$.
Граф порядка $n$ назовем \emph{свитчингово разделимым},
если существует отделимое множество вершин.
Далее слово ``свитчингово'' будем опускать.

\begin{remark}
Если граф разделим, то все его свитчинги, а также его дополнение, разделимы. Все графы порядка $4$ разделимы.
\end{remark}

Целью настоящей работы является изучение взаимоотношений между разделимостью графа и разделимостью его подграфов. Мотивацией исследования является связь с разделимостью $n$-арных квазигрупп и булевых функций, описанная в разделе~\ref{s:gbq}. В разделах~\ref{s:n-3} и~\ref{s:n-2} мы  докажем следующие две теоремы:
%:th:n-3
\begin{theorem}\label{th:n-3}
Если все подграфы порядка $n-1$ и $n-2$ графа $G=(V,E)$ порядка $n$ разделимы, то $G$ --- разделимый граф.
\end{theorem}
%:th:n-2
\begin{theorem}\label{th:n-2}
Для любого нечетного $n$ существует неразделимый граф порядка $n$, у которого все подграфы порядка $n-1$ разделимы.
\end{theorem}
Вопрос с четным порядком остается открытым, однако полным перебором для порядков $6$ и $8$ таких графов не найдено.
\begin{conjecture}
Из любого неразделимого графа четного порядка удалением одной вершины можно получить неразделимый подграф.
\end{conjecture}
Основной результат можно переформулировать следующим образом:
\begin{corollary}
Из любого неразделимого графа удалением одной или двух вершин можно получить неразделимый подграф, причем удаления одной вершины не всегда достаточно.
\end{corollary}
%:222222222222222222222222222222222222222222222222222222222
\section{Доказательство теоремы \ref{th:n-3}}\label{s:n-3}
Пусть $\kappa $ --- максимальный порядок неразделимого собственного подграфа графа $G$, и $K$ --- множество вершин некоторого такого подграфа.
По условию $3\leq \kappa\leq n-3$

Сначала рассмотрим случай $\kappa>3$. Это неравенство будет неявно спользоваться в местах, где выводится противоречие с неразделимостью $K$.
Для четырех вершин $a$, $b$, $c$, $d$ графа $G$ через $N(a,b;c,d)$ обозначим число ребер
графа $G$ среди $\{a,c\}$, $\{a,d\}$, $\{b,c\}$, $\{b,d\}$.
%:L:even
\begin{lemma}\label{l:even}
Множество вершин $W$ графа $G$ отделимо тогда и только тогда, когда
для любых попарно различных $a,b$ из $W$ и $c,d$ из $V\setminus W$ число $N(a,b;c,d)$ четно.
\end{lemma}
\begin{proof}
{\em Только тогда.} Пусть $W$ отделимо, причем множество вершин $U$ задает отделяющий свитчинг. Для любых попарно различных $a,b$ из $W$ и $c,d$ из $V\setminus W$ среди ребер $\{a,c\}$, $\{a,d\}$, $\{b,c\}$, $\{b,d\}$ графу $G$ принадлежат в точности те же ребра, что и полному двудольному графу $K_{U,V\setminus U}$. Легко убедиться, что в полном двудольном графе число ребер, соединяющих две пары вершин, всегда четно.

{\em Тогда.} Рассмотрим две несмежные вершины $a$ из $W$ и $c$ из $V\setminus W$ (если таких нет, то $W$-свитчинг отделяет множество $W$). Введем четыре множества
$$ W_0 = \{ b \in W \mid \{b,c\}    \in E \}, \qquad
   V_0 = \{ d \in V\setminus W \mid \{a,d\}    \in E \}, $$
$$ W_1 = \{ b \in W \mid \{b,c\}\not\in E \}, \qquad
   V_1 = \{ d \in V\setminus W \mid \{a,d\}\not\in E \}. $$
Мы утверждаем, что $b$ из $W_i$ и $d$ из $V_j$ соединены ребром если и только если $i+j=1$. (Действительно, если $b=a$ или $d=c$, то это следует прямо из определений множеств $W_i$ и $V_j$; в противном случае --- из четности $N(a,b;c,d)$.) Таким образом, взяв $U=W_0\cup V_0$, мы получаем $U$-свитчинг без ребер, соединяющих $W$ и $V\setminus W$. Лемма~\ref{l:even} доказана.\end{proof}
%:L:4+4
\begin{lemma}\label{l:4+4}
Для любых попарно различных вершин $a$, $b$, $c$, $d$, $e$ графа $G$
из четности $N(a,b;c,d)$ и $N(a,b;c,e)$ следует четность $N(a,b;d,e)$.
\end{lemma}
\begin{proof}
$N(a,b;d,e)=N(a,b;c,d)+N(a,b;c,e)-2|\{\{a,c\},\{b,c\}\}\cap E|$.
\end{proof}
Рассмотрим произвольную вершину $v$ не из $K$. По определению $K$ граф $G|_{K\cup \{v\}}$ разделим, то есть $v$ принадлежит некоторому отделимому в $G|_{K\cup \{v\}}$ множеству вершин. Если это множество имеет больше двух вершин, то граф $G|_{K}$ также разделим, что противоречит определению $K$. Таким образом, для любой вершины $v\not\in K$ найдется вершина $u=u(v)\in K$ такая, что $\{v,u\}$ отделимо в $G|_{K\cup \{v\}}$. Причем $u(v)$ определена однозначно (если $\{v,u\}$ и $\{v,u'\}$ отделимы и $u\neq u'$, то из лемм~\ref{l:even} и~\ref{l:4+4} получаем отделимость $\{v,u,u'\}$ и противоречие с неразделимостью $G|_{K}$).
%:Prop:vv'
\begin{proposition}\label{p:1}
Для любых вершин $v,v'\not\in K$ таких, что $u(v)\neq u(v')$,
число $N(v,u(v);v',u(v'))$ четно.
\end{proposition}
\begin{proof}
Из условия следует, что подграф $G|_{K\cup \{v,v'\}}$ разделим.
Пусть $v$ принадлежит множеству вершин $M$, отделимому в этом подграфе.
Рассмотрим подслучаи.

1) $|M\cap K|=0$, т.\,е. $M=\{v,v'\}$.
По лемме~\ref{l:even} для любых $c,d\in K\setminus \{u(v)\}$ число $N(v,v';c,d)$ четно. Но то же верно для $N(v,u(v);c,d)$, поскольку
$\{v,u(v)\}$ отделимо в $G|_{K\cup \{v\}}$. И по лемме~\ref{l:4+4}
число $N(v',u(v);c,d)$ также четно, откуда следует
отделимость $\{v,v',u(v)\}$ в $G|_{K\cup \{v,v'\}}$,
отделимость $\{v',u(v)\}$ в $G|_{K\cup \{v'\}}$
и противоречие с однозначностью $u(v')$.

2) $|M\cap K|=1$. Если $M=\{v,u(v)\}$, то утверждение предложения~\ref{p:1} следует из леммы~\ref{l:even}. Любой другой случай противоречит однозначности $u(v)$ или $u(v')$.

3) $2\leq |M\cap K|\leq \kappa-2$ противоречит неразделимости $G|_{K}$.

4) $|M\cap K|=\kappa-1$. Единственной вершиной из $K\setminus M$ может быть только $u(v')$, и утверждение предложения~\ref{p:1} следует из леммы~\ref{l:even}.

Предложение~\ref{p:1} доказано.
\end{proof}
Рассмотрим некоторую вершину $w$ из $K$ с непустым прообразом $u^{-1}(w)$.
Обозначим $W = u^{-1}(w)\cup\{w\}$ и покажем, что это множество отделимо. По лемме~\ref{l:even} это эквивалентно четности $N(a,b;c,d)$ для любых попарно различных вершин $a,b$ из $W$ и $c,d$ из $V\setminus W$. По лемме~\ref{l:4+4} достаточно рассмотреть случай $b=w=u(a)$.
Если $c$ и $d$ принадлежат $K$, то нужная четность следует из определения $u(a)$; если $d=u(c)$ --- то из предложения~\ref{p:1}.
Остальные случаи выводятся из этих двух посредством леммы~\ref{l:4+4}. Таким образом, разделимость $G$ в случае $\kappa>3$ доказана.

Рассмотрим случай $\kappa=3$.
Без потери общности можно считать, что граф $G$ содержит некоторую изолированную вершину $o$ (в противном случае выберем произвольную вершину в качестве $o$ и рассмотрим свитчинг по множеству смежных с ней вершин).
%:Prop:avoid
\begin{proposition}\label{p:14}
Граф $G$ не содержит подграфа вида
$$\left(\{o,a,b,c,d\},\left\{\{a,b\},\{b,c\},\{c,d\}\right\}\right).$$
\end{proposition}
\begin{proof}
Существование такого подграфа противоречило бы условию теоремы в силу его неразделимости.
\end{proof}

Две несовпадающие вершины назовем {\em парными}, если любая вершина, смежная с одной из них, смежна и с другой.
%:Prop:twin
\begin{proposition}\label{p:twin}
Если $v$, $w$ --- парные вершины, то $\{v,w\}$ --- отделимое множество.
\end{proposition}
\begin{proof}
Пусть $U$ --- множество вершин, смежных с $v$ и $w$. Тогда $U$-свитчинг графа $G$ не содержит ребер, соединяющих $v$ или $w$ с остальными вершинами.
\end{proof}

Таким образом, чтобы показать разделимость графа $G$, достаточно найти две парные вершины.
Рассмотрим максимальную последовательность вершин\linebreak
$\bar u = (u_1, u_2, \ldots, u_t)$,
удовлетворяющую следующему свойству:

(*) вершины $u_i$ и $u_j$, $1\leq i<j\leq t$, смежны если и только если $i$ нечетно.

Если $t=1$, то граф пустой и доказывать нечего.
Пусть $t>1$.
Покажем от противного, что вершины $u_{t-1}$ и $u_t$ парные. Предположим, что это не так, тогда найдется вершина $w$, смежная ровно с одной из $u_{t-1}$, $u_t$, причем из (*) следует, что $w$ не принадлежит $\bar u$. Мы покажем, что в этом случае последовательность $\bar u$ не максимальна.

Сначала рассмотрим случай, когда $t$ нечетно. Без потери общности можно считать, что $w$ смежна с $u_t$ (в противном случае мы можем переставить
местами $u_{t-1}$ и $u_t$, при этом свойство (*) сохранится). Заметим следующее:

\begin{itemize}
\item[--]
 для каждого нечетного $i$ меньше $t$ вершины $u_i$ и $w$ смежны, иначе вершины $o$, $u_{t-1}$, $u_i$, $u_t$, $w$ порождают запрещенный подграф (предложение~\ref{p:14});
\item[--]
 для каждого четного $i$ меньше $t-1$ вершины $u_i$ и $w$ несмежны, иначе вершины $o$, $u_i$, $w$, $u_{t-2}$, $u_{t-1}$ порождают запрещенный подграф (предложение~\ref{p:14}).
\end{itemize}

Таким образом, последовательность $(u_1, u_2, \ldots, u_t,w)$
опровергает максимальность $\bar u$ среди последовательностей, удовлетворяющих свойству (*).
Полученное противоречие доказывает, что вершины $u_{t-1}$ и $u_t$ парные.

Случай четного $t$ рассматривается аналогично.
Теорема доказана.
%:333333333333333333333333333333333333333333333333333333333
\section{Доказательство теоремы \ref{th:n-2}}\label{s:n-2}
Пусть $n$ нечетно. Обозначим через $G_n$ граф с множеством вершин $V_n=\{v_i\}_{i=0}^{n-1}$ и ребрами
$\{v_i,v_{i+j}\}$, $j=1,\ldots,\left\lfloor  {n \over 4}\right\rfloor$
(здесь и далее вычисления с индексами будем производить по модулю $n$).
%:Prop:nonseparable
\begin{proposition}\label{p:nonsep}
Граф $G_n$ неразделим.
\end{proposition}
\begin{proof}
Обозначим $m= \lfloor (n+1)/4 \rfloor $ и $u_i=v_{im}$.
Поскольку $n=4m\pm
 1$, числа $m$ и $n$ взаимно просты, откуда
$\{u_i\}_{i=0}^{n-1} = V_n$. Теперь рассмотрим произвольное подмножество
$A\subset V_n$ мощности не меньше $2$ и не больше $n-2$ и покажем, что оно неотделимо (эквивалентно, $V_n\setminus A$ неотделимо).
Имеет место один из следующих двух случаев:

1) Для некоторого $i$
либо $u_{i},u_{i+1}\in A$, $u_{i+2},u_{i+3}\not\in A$,
либо $u_{i},u_{i+1}\not\in A$, $u_{i+2},u_{i+3}\in A$.
Тогда $N(u_{i},u_{i+1},u_{i+2},u_{i+3})=1$, см. рис.~\ref{f:}, и по лемме~\ref{l:even} множество $A$ неотделимо.
%:Figure
\begin{figure}
\begin{center}
\def\iiii(#1-#2-#3-#4){\path[every node/.style={draw,rectangle,inner sep=1.2pt,fill=white}]
(#1) node{$\ a_i\ $}
(#2) node{$a_{\scriptscriptstyle i{+}1}$}
(#3) node{$a_{\scriptscriptstyle i{+}2}$}
(#4) node{$a_{\scriptscriptstyle i{+}3}$};}
\newcommand\dvsr{3}
\def\DefNode #1(#2){\node(#1) at ($(10,10)+sin(#2*360/\dvsr)*(10,0)+cos(#2*360/\dvsr)*(0,10)$){};}
\def\dfNM#1#2{\renewcommand\dvsr{#1}\node[rectangle,draw=white] at (10,10) {$\displaystyle{ n=#1 \atop m=#2}$};}
%55555555555555555555555555555555555555555555555555555555555
\mbox{}\hfill
\begin{tikzpicture}[scale=0.24, every node/.style={circle,draw,inner sep=2.0pt}]
\dfNM{13}{3}
\foreach \nam / \num in {a/0,b/1,c/2,d/3,e/4,f/5,g/6,h/7,i/8,j/9,k/10,l/11,m/12} {\DefNode \nam(\num)}
\draw[line width=0.7pt]
 (a)--(b)--(c)--(d)--(e)--(f)--(g)--(h)--(i)--(j)--(k)--(l)--(m)--(a)
    --(c)--(e)--(g)--(i)--(k)--(m)--(b)--(d)--(f)--(h)--(j)--(l)--(a)
    --(d)--(g)--(j)--(m)--(c) (l)--(b)--(e)--(h)--(k)--(a);
\draw[line width=2.2pt] (c)--(f)--(i)--(l);
\iiii(c-f-i-l)
\end{tikzpicture}
%66666666666666666666666666666666666666666666666666666666666
\ \ \ \ \ \ \ \ %
\begin{tikzpicture}[scale=0.24, every node/.style={circle,draw,inner sep=2.0pt}]
\dfNM{15}{4}
\foreach \nam / \num in {a/0,b/1,c/2,d/3,e/4,f/5,g/6,h/7,i/8,j/9,k/10,l/11,m/12,n/13,o/14} {\DefNode \nam(\num)}
\draw[line width=0.7pt]
 (a)--(b)--(c)--(d)--(e)--(f)--(g)--(h)--(i)--(j)--(k)--(l)--(m)--(n)--(o)--(a)
    --(c)--(e)--(g)--(i)--(k)--(m)--(o)--(b)--(d)--(f)--(h)--(j)--(l)--(n)--(a)
(a)--(d)--(g)--(j)--(m)--(a)(c)--(f)--(i)--(l)--(o)(b)--(e)--(h)--(k)--(n)--(b);
\draw[line width=2.2pt] (o)--(c);
\iiii(c-g-k-o)
\end{tikzpicture}
\hfill\mbox{}
 \caption{\label{f:} 	Примеры графа $G_n$, случаи $n\equiv 1\bmod 4$ и $n\equiv 3\bmod 4$.
  }
\end{center}
\end{figure}
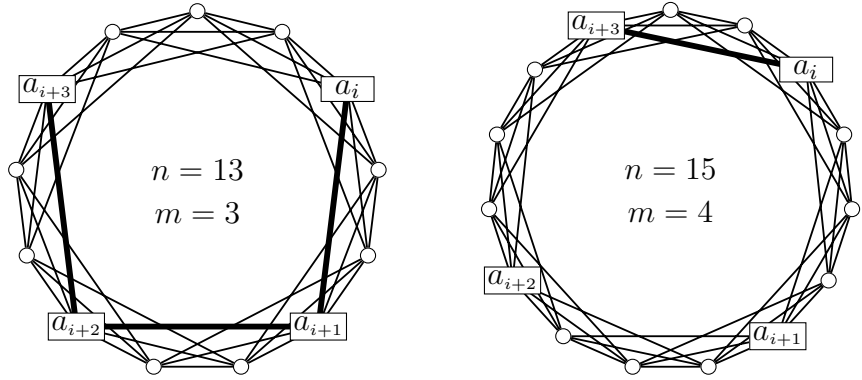

2) Для некоторого $i$
либо $u_{i},u_{i+2}\in A$, $u_{i+1}\not\in A$,
либо $u_{i},u_{i+2}\not\in A$, $u_{i+1}\in A$. Рассмотрим для примера второй подслучай. Заметим, что любая вершина $u_j$, отличная от $u_i$, $u_{i+1}$, $u_{i+2}$, смежна ровно с одной из $u_i$, $u_{i+2}$ (см. рис.~\ref{f:}). Взяв такую $u_j$ из $A$, получаем, что
$N(u_{i},u_{i+2};u_{i+1},u_{j})=1$ или $3$, в зависимости от $n\bmod 4$ (см. рис.~\ref{f:}), и по лемме~\ref{l:even} множество $A$ неотделимо.

По определению если любое множество вершин неотделимо, то граф является неразделимым. Предложение~\ref{p:nonsep} доказано.
\end{proof}
%:Prop:separable
\begin{proposition}\label{p:sep}
Удаление любой вершины в графе $G_n$ приводит к разделимому графу.
\end{proposition}
\begin{proof}
В силу симметрии можно считать, что удалили вершину $v_0$.
Легко видеть (рис.~\ref{f:}), что все оставшиеся вершины кроме $v_m$ и $v_{-m}$ делятся
на смежные с $v_m$ и смежные с $v_{-m}$. Откуда следует, множество
$\{v_m,v_{-m}\}$ отделимо (для соответствующего $U$-свитчинга нужно взять множество $U$, состоящее из $v_m$ и всех вершин, смежных с $v_{-m}$).
\end{proof}
Таким образом, теорема~\ref{th:n-2} доказана. Заметим, что в силу теоремы~\ref{th:n-3} в графе $G_n$ есть неразделимый подграф порядка $n-2$.
Аналогично доказательству предложения~\ref{p:nonsep} можно показать, что удаление вершин $v_i$ и $v_{i+m}$ приводит к неразделимому подграфу.
%:444444444444444444444444444444444444444444444444444444444
\section{Графы, булевы функции, квазигруппы}\label{s:gbq}
В этом разделе мы кратко обсудим связь разделимости графов с аналогичным свойством для булевых функций и $n$-арных квазигрупп.
Подграфам графа соответствуют так называемые ретракты $n$-арных квазигрупп и подфункции булевых функций, и то и то получается фиксацией некоторых аргументов. В терминах ретрактов и подфункций для $n$-арных  квазигрупп и булевых функций верны теоремы, аналогичные теоремам~\ref{th:n-3} и ~\ref{th:n-2} (последняя известна для $n$-арных  квазигрупп только если порядок кратен $4$). Причем, учитывая, что по графу при помощи квадратичного многочлена можно построить булеву функцию и затем $n$-арную квазигруппу порядка $4$, теорема~\ref{th:n-3} является, вообще говоря, следствием соответствующей теоремы для квазигрупп \cite{Kro:n-3,KroPot:order4}, а из теоремы~\ref{th:n-2}, наоборот, следует существование аналогичного примера в $n$-арных  квазигруппах порядка $4$ \cite{Kro:n-2}.

\subsection{Расширенные булевы функции}
\emph{Расширенной булевой функцией} назовем частичную булеву функцию, заданную на наборах с четным числом единиц. (Заметим, что расширенную булеву функцию можно интерпретировать как обычную булеву функцию от на единицу меньшего числа аргументов.) Расширенную булеву функцию $f$ от $n$ аргументов назовем \emph{разделимой}, если она представима в виде суммы двух булевых функций $f'$ и $f''$ от $n-2$ или меньшего числа аргументов, причем наборы аргументов $f'$ и $f''$ не перекрываются (ограничение $n-2$ достаточно естественно: в этом случае булевы функции в разложении могут быть заданы меньшим числом значений в точках, чем сама расширенная булева функция). \emph{Степенью} расширенной булевой функции назовем минимальную степень многочлена (над полем GF(2)), с помощью которого она может быть представлена. Под термином <<квадратичный>> будем подразумевать <<степени не больше двух>>.
Графом квадратичного многочлена назовем граф на множестве аргументов, у которого две вершины смежны тогда и только тогда, когда произведение соответствующих переменных входит в многочлен.
%:Lem:bol-swe
\begin{lemma}\label{l:bol-swe}
Множество графов, соответствующих представлениям данной квадратичной расширенной булевой функции в виде квадратичного полинома, образует свитчинговый класс.
\end{lemma}
\begin{proof}
Любая (и в частности, квадратичная) расширенная булева функция $f$ от $n$ аргументов, будучи булевой функцией от первых $n-1$ своих аргументов, единственным образом представима в виде
$$ f(x_1,\ldots,x_{n-1},x_{n})=p(x_1,\ldots,x_{n-1}) ,$$
где $p$ --- многочлен.

Любой  многочлен $r$ от $n$ переменных $x_1,\ldots,x_{n-1},x_{n}$
однозначно представим в виде
$$ q(x_1,\ldots,x_{n-1})+(x_1+\dots+x_{n-1}+x_{n})l(x_1,\ldots,x_{n-1}),$$
где $q$ и $l$ --- многочлены от $x_1,\ldots,x_{n-1}$, причем если многочлен $r$ квадратичный, то $q$ квадратичный и $l$ линейный.
Поскольку $x_1+\dots+x_{n-1}+x_{n}=0$ везде на области определения расширенной булевой функции,
многочлен $q$ совпадает у всех многочленов, представляющих одну и ту же расширенную булеву функцию.
Легко убедиться, что добавление $(x_1+\dots+x_{n-1}+x_{n})l(x_1,\ldots,x_{n-1})$ с линейным $l$ приводит к свитчингу соответствующего графа, точнее, к $U$-свитчингу, где $U$ --- множество переменных, от которых $l$ существенно зависит. Отсюда следует утверждение леммы.
\end{proof}
\begin{lemma}\label{l:bol-sep}
Квадратичная расширенная булева функция разделима тогда и только тогда,
когда разделимы графы квадратичных многочленов, представляющих эту функцию.
\end{lemma}
\begin{proof}
Из разделимости графа по определению следует разделимость булевой функции, представимой квадратичным многочленом с соответствующим графом.

Для доказательства обратного, с учетом предыдущей леммы, достаточно показать,
что для разделимой квадратичной расширенной булевой функции $f$ элементы ее некоторого разложения $f'$ и $f''$ из определения разделимости также квадратичны.
Пусть $$f(x_1,\ldots,x_n)=f'(\bar y)+f''(\bar z),$$
где $\bar y$ и $\bar z$ -- непересекающиеся наборы переменных из $x_1,\ldots,x_n$.
Представим $f'(\bar y)+f''(\bar z)$, как везде определенную булеву функцию, в виде
$$f'(\bar y)+f''(\bar z)=q(x_1,\ldots,x_{n-1})+(x_1+\dots+x_{n})l(x_1,\ldots,x_{n-1}).$$
Поскольку $f$ квадратична, $q$ также квадратична (см. доказательство предыдущей леммы).
Разобьем полином $l$ в сумму двух $l_1+l_2$, где $l_1$ линейный, а $l_2$ составлен из мономов степени $2$ и выше. Имеем
$$f'(\bar y)+f''(\bar z)=q(x_1,\ldots,x_{n-1})+\sum_{i=1}^{n}x_{i}l_1(x_1,\ldots,x_{n-1})+\sum_{i=1}^{n}x_{i}l_2(x_1,\ldots,x_{n-1}).$$
Легко видеть, что последнее слагаемое состоит из мономов степени не меньше $3$ (действительно, если $l_2$ содержит моном $x_i x_j$, то в произведении с $x_i+x_j$ он даст нуль, а произведение с остальными переменными даст мономы третьей степени).
Таким образом, отбрасывание этого слагаемого равносильно отбрасыванию мономов степени больше двух в полиномиальном представлении $f'$ и $f''$, после чего имеем
$$g'(\bar y)+g''(\bar z)=q(x_1,\ldots,x_{n-1})+\sum_{i=1}^{n}x_{i}l_1(x_1,\ldots,x_{n-1}),$$
для некоторых квадратичных функций $g'$ и $g''$. Очевидно, что $g'(\bar y)+g''(\bar z)$ также равны
расширенной булевой функции $f$ на всей ее области определения, то есть мы получили квадратичное представление $f$, которому соответствует разделимый граф. По лемме~\ref{l:bol-swe} все другие квадратичные представления также соответствуют разделимым графам.
\end{proof}
\subsection{$n$-Арные квазигруппы}
Пусть $\Sigma$ --- некоторое множество. $n$-Арная операция $Q:\Sigma^n \to \Sigma$ называется
\emph{$n$-арной квазигруппой} порядка $|\Sigma|$, если в уравнении $x_0=Q(x_1,\ldots,x_n)$ значения любых $n$ переменных однозначно задают значение оставшейся переменной.
(Строго говоря, $n$-арной квазигруппой называется пара $(\Sigma,Q)$, наше определение --- общепринятое упрощение терминологии.)
Из определения следует, что $n$-арная квазигруппа обратима в каждой позиции,
в случае конечного порядка это свойство можно взять за определение.
Введем обозначение $Q\langle x_0,x_1,\ldots,x_n\rangle \Leftrightarrow x_0=Q(x_1,\ldots,x_n)$
для предикатной записи $n$-арной квазигруппы; часто предикатная запись удобнее функциональной ввиду симметричности относительно всех переменных.
Если в предикате $Q\langle \ldots\rangle$ зафиксировать значения некоторых $m\in\{1,\ldots,n\}$ аргументов, то полученный $(n{+}1{-}m)$-местный предикат соответствует некоторой $(n{-}m)$-арной квазигруппе, которая называется \emph{ретрактом} квазигруппы $Q$.
$n$-Арная квазигруппа называется \emph{разделимой}, если она представима в виде безповторной суперпозиции двух квазигрупп меньшей арности, где порядок переменных в суперпозиции может отличаться от первоначального.

\begin{remark}
  В литературе также известен термин <<приводимая $n$-арная квазигруппа>>, который чаще относится к представимости в виде суперпозиции с тем же порядком переменных. В англоязычной литературе перевод слова <<разделимая>> сильно перегружен, поэтому разделимые квазигруппы также называют <<перестановочно приводимыми>>.
\end{remark}

Пусть $\Sigma=\{[0,0],[0,1],[1,0],[1,1]\}$ --- множество двоичных пар и $\lambda$ --- некоторая расширенная булева функция от $n+1$ переменной. Предикат
$$
Q_{\lambda}\Bigl\langle[x_0,y_0],\ldots,[x_n,y_n]\Bigr\rangle
\Leftrightarrow
\left\{{
|x_{0}+\dots+x_{n}| =0,\hfill\mbox{}
\atop
|y_{0}+\dots+y_{n}| = \lambda(x_{0},\dots,x_{n})
}\right.
$$
соответствует $n$-арной квазигруппе $Q$ (сложение производится по модулю $2$), эта конструкция является частным случаем сплетения $n$-арных квазигрупп, в данном случае тривиальных квазигрупп порядка $2$. (Термин <<сплетение>> для $n$-арных квазигрупп не соответствует сплетению групп, поэтому следует быть осторожным при его использовании, дабы избежать возможных разночтений.)

\begin{lemma}[\cite{Kro:n-2}]\label{l:bol-qua}
  Разделимость $n$-арной квазигруппы $Q_{\lambda}$ эквивалентна разделимости расширенной булевой функции $\lambda$.
\end{lemma}

Таким образом, разделимость графов тесно связана с разделимостью $n$-арных квазигрупп, по крайней мере в рамках следующей конструкции: по графу порядка $n+1$ мы строим квадратичную расширенную булеву функцию $\lambda$ (ребру соответствует моном степени $2$, линейная часть выбирается произвольно); потом строим $n$-арную квазигруппу $Q_\lambda$ порядка $4$; далее можем построить $n$-арную квазигруппу порядка $4k$ для любого, в том числе бесконечного, $k$, при помощи прямого произведения с $n$-арной квазигруппой $P(x_1,\ldots,x_n)=x_1\star\ldots\star x_n$, где $\star$ --- коммутативная групповая операция.

При этом разделимость каждого звена цепочки эквивалентна разделимости всех остальных звеньев. Более того, разделимость подграфа эквивалентна разделимости соответствующих подфункций расширенной булевой функции и ретрактов $n$-арной квазигруппы. Отсюда теорема~\ref{th:n-3} есть следствие приведенной ниже теоремы~\ref{th:qua-n-3}, а приведенная ниже  теорема~\ref{th:qua-n-2} следует из теоремы~\ref{th:n-2}. В защиту разделов~\ref{s:n-3} и~\ref{s:n-2} стоит заметить, что содержащиеся в них доказательства значительно проще имеющихся доказательств для квазигрупп.

В заключение сформулируем известные теоремы для $n$-арных квазигрупп, связанные с рассмотренными в настоящей работе теоремами~\ref{th:n-3} и~\ref{th:n-2}. Для $n$-арной квазигруппы $Q$ обозначим через $\kappa(Q)$ наибольшую арность ее неразделимого ретракта.

\begin{theorem}\label{th:qua-n-3}
Если $\kappa(Q)<n-2$, то $n$-арная квазигруппа $Q$ разделима.
Если $\kappa(Q)=n-2>2$ и порядок $Q$ --- простое число, то $Q$ разделима.
\end{theorem}
Случай $2<\kappa(Q)<n-2$ доказан в \cite{Kro:n-3}, случай $\kappa(Q)=2$ для порядка $4$ --- в \cite{KroPot:order4}, общий случай $\kappa(Q)=2$ и случай $\kappa(Q)=n-2$ для простого порядка --- в \cite{KroPot:ir2}. Теорема~\ref{th:qua-n-3} полезна при индуктивной характеризации классов $n$-арных квазигрупп, например, при доказательстве, что любая $n$-арная квазигруппа порядка $4$ полулинейна (то есть эквивалентна некоторой $Q_\lambda$) или разделима использовалось, что минимальный потенциальный контрпример обязан иметь неразделимый полулинейный $(n{-}1)$- или $(n{-}2)$-арный ретракт.
\begin{theorem}\label{th:qua-n-2}
Для любого четного $n$ и любого $k$ существует неразделимая $n$-арная квазигруппа $Q$ порядка $4k$ такая, что  $\kappa(Q)=n-2$ {\rm \cite{Kro:n-2}}. Для любого $n\geq 3$ и $k\geq 4$ существует неразделимая $n$-арная квазигруппа порядка $k$ с $\kappa(Q)=n-1$ {\rm \cite{KPS}}.
\end{theorem}
Теоретически на предмет существования неразделимой $n$-арной квазигруппы остались неисследованы следующие случаи:

$\kappa(Q)=n-2$, нечетное $n$ (связан с гипотезой, сформулированной во введении), непростой порядок (для простого --- не существует);

$\kappa(Q)=n-2$, произвольное $n$, непростой порядок некратный $4$;

$\kappa(Q)=n-2=2$, $n=4$, порядок некратный $4$ (для кратного --- существует), в частности простой порядок.

Известен пример неразделимой $4$-арной квазигруппы $Q$ порядка $6$ с $\kappa(Q)=2$.
\vskip 5mm

Автор благодарит А.\,Н.\,Глебова, В.\,Н.\,Потапова, А.\,В.\,Пяткина и анонимного рецензента за интерес к данной работе и замечания, благодаря которым была обнаружена неполнота доказательства теоремы~\ref{th:n-3} в первом варианте манускрипта.